\documentclass[12pt]{article}

\usepackage{amsmath,amssymb}

\newtheorem{thm}{Theorem}
\newtheorem{prop}{Proposition}
\newtheorem{lem}{Lemma}

\newtheorem{cor}{Corollary}

\title{A characterization\\ of Gorenstein toric Del Pezzo $n$-folds
}

\author{Shoetsu Ogata\thanks{e-mail:  ogata{\char'100}math.tohoku.ac.jp} and  
Huai-Liang Zhao\thanks{e-mail: sa9m37{\char'100}math.tohoku.ac.jp}\\
Mathematical Institute, Tohoku University\\ Sendai 980-8578, Japan}

\date{}

\begin{document}
\maketitle


\begin{abstract}
We give a characterization of Gorenstein toric Fano $n$-fold with index $n-1$, which is called
Gorenstein toric Del Pezzo $n$-folds, among toric varieties.
In practice, we obtain a condition for a lattice $n$-polytope to be a Gorenstein Fano polytope.
In our proof we do not use the Batyrev-Juny's classification of Gorenstein Del Pezzo
$n$-polytopes.
\end{abstract}
\section*{Introduction}

A nonsingular projective variety $X$ is called {\it Fano} if its anti-canonical divisor
$-K_X$ is ample, and the number $i_X:=\max\{i\in \mathbb{N}; -K_X=iD \ \mbox{for a divisor
$D$}\}$ is called the {\it Fano index}, or simply, {\it index} of $X$.
Even if a variety $X$ has at worst  Gorenstein singularities, 
we can define it to be Fano and its index.

A nonsingular Fano variety with index $i_X=\dim X-1$ is called a {\it Del Pezzo} manifold.
Fujita \cite{Fj1} \cite{Fj2} classifies Del Pezzo manifolds.
Batyrev and Juny \cite{BJ} classified Gorenstein toric Del Pezzo varieties.

In this paper, we  give a characterization of Gorenstein toric Del Pezzo varieties
among toric varieties.

\begin{thm}\label{0;t0}
Let $X$ be a projective toric surface.  If $X$ has an ample line bundle $L$ with
$\dim \Gamma(X, L\otimes \omega_X)=1$, then it is a Gorenstein toric Fano surface,
where $\omega_X$ is the dialyzing sheaf of $X$.
\end{thm}
\begin{thm}\label{0;t1}
Let $X$ be a projective toric variety of dimension $n$ with $n\ge3$.
If $X$ has an ample line bundle $L$  with $\dim \Gamma(X,
L^{\otimes (n-1)}\otimes \omega_X)=1$,
then $X$ is a Gorenstein toric Del Pezzo variety.
\end{thm}

We note that the condition in Theorem~\ref{0;t1} is not trivial 
for $X$ to be Gerenstein. 
Ogata and Zhao \cite{OZ} give examples of
 polarized toric $n$-fold $(X, L)$ with $\dim \Gamma(X, L^{\otimes (n-2)}\otimes 
\omega_X)=1$ which is not Gorenstein for every $n$ with $n\ge4$.  

Ogata and Zhao \cite{OZ} give a characterization of a  toric variety $X$
 to be Gorenstein toric Fano with index $i_X=\dim X$.
They
 also give a characterization of a special class of Gorenstein toric Del Pezzo varieties.

\begin{thm}[Ogata-Zhao]\label{0;t01}
Let $X$ be a projective toric variety of dimension $n$ with $n\ge3$.
If $X$ fas an ample line bundle $L$ satisfying the condition that $\dim \Gamma(X, L)=
n+1$ and $\dim \Gamma(X, L^{\otimes(n-1)}\otimes \omega_X)=1$, then
it is a Gorenstein toric Del Pezzo variety.
\end{thm}

In \cite{OZ} they remarked that the line bundle $L$ in Theorem~\ref{0;t01} is not very ample.
In this paper we show the normal generation for $L$ with more global sections.

\begin{thm}\label{0;t02}
Let $X$ be a projective toric variety of dimension $n$ with $n\ge3$.
If  an ample line bundle $L$ on $X$ satisfies the condition that $\dim \Gamma(X, L)\ge
n+2$ and $\dim \Gamma(X, L^{\otimes(n-1)}\otimes \omega_X)=1$, then
$L$ is normally generated.
\end{thm}
This Theorem is a corollary of Proposition~\ref{1;p1}, which will be proved
in terms of algebraic geometry in Section 1.

All statements in the above Theorems~\ref{0;t0}, \ref{0;t1} and \ref{0;t01}
 can be interpreted in terms of lattice polytopes.
In the next section, we recall the fundamental notions of lattice polytopes and
the relationship between polarized toric varieties and lattice polytopes.

Theorem~\ref{0;t0} is given as Proposition~\ref{2;p1} and Theorem~\ref{0;t1} is separated 
into the case of dimension three as Proposition~\ref{3;p1}
and the case of higher dimension as Proposition~\ref{4;p1}.
In our proof, we do not use the classification of Gorenstein Del Pezzo $n$-polytopes
given by Batyrev and Juny \cite{BJ}.

\section{Toric Varieties and Lattice Polytopes}

Let $M=\mathbb{Z}^n$ be a free abelian group of rank $n$ and $M_{\mathbb{R}}:=
M\otimes_{\mathbb{Z}}\mathbb{R} \cong\mathbb{R}^n$ the extension of coefficients
into real numbers.
We define a {\it lattice polytope} $P$  in $M_{\mathbb{R}}$ as the convex hull 
$P:= \mbox{Conv}\{m_1,
\dots, m_r\}$ of a finite subset $\{m_1, \dots, m_r\}$ of $M$.
We define the dimension of a lattice polytope $P$ as that of the smallest
affine subspace containing $P$.  A lattice polytope of dimension $n$ is 
sometimes called a lattice $n$-polytope.
In particular, a lattice $n$-simplex with only $n+1$ vertices is 
called {\it lattice $n$-simplex}.

The space $\Gamma(X, L)$ of global sections of an 
ample line bundle $L$ on a toric variety $X$ of dimension $n$
 is parametrized by the set of lattice points in a lattice polytope $P$
of dimension $n$ (see, for instance, Oda's book\cite[Section 2.2]{Od} 
or Fulton's book\cite[Section 3.5]{Fu}). 
 And $k$ times tensor product $L^{\otimes k}$ corresponds to the polytope
 $kP:=\{kx\in M_{\mathbb{R}}; x\in P\}$.   Furthermore, the surjectivity of
 the multiplication map $\Gamma(X, L^{\otimes k})\otimes \Gamma(X, L) \to
 \Gamma(X, L^{\otimes(k+1)})$
 is equivalent to the equality
 \begin{equation}\label{eq1}
 (kP)\cap M +P\cap M = ((k+1)P)\cap M.
 \end{equation}
 Thus, $L$ is normally generated if and only if the equality (\ref{eq1}) holds for all
 $k\ge1$.  In that case, we define as $P$ is {\it normal}.
In particular, the dimension of the space of global sections of $L\otimes \omega_X$ is
equal to the number of lattice points contained in the interior of $P$, i.e.,
\begin{equation}\label{eq2}
\dim \Gamma(X, L\otimes \omega_X)=\sharp\{(\mbox{Int}P)\cap M\}.
\end{equation}

Here we can prove Theorem~\ref{0;t02}, which is the case $g=1$ in the following proposition.

\begin{prop}\label{1;p1}
Let $L$ be an ample line bundle on a projective toric variety $X$ of dimension $n$ with $n\ge3$.
Assume that $\Gamma(X, L^{\otimes(n-2)}\otimes \omega_X)=0$.  Set
$g:=\dim \Gamma(X, L^{\otimes(n-1)}\otimes \omega_X)$.
If $\dim \Gamma(X, L)\ge n+g+1$, then $L$ is normally generated. 
\end{prop}
{\it Proof}.   Let $H_1, \dots, H_{n-1}$ be general members of the linear system $|L|$ and
$C:=H_1\cap \dots \cap H_{n-1}$ the intersection curve.
Since $X$ is a normal variety, $C$ is nonsingular.  Let $L_C$ be the restriction of $L$ to $C$.

Since $\chi(X, L^{\otimes r}\otimes \omega_X)=0$ for $1\le r\le n-2$, we have
$\chi(C, L_C)=\chi(X, L) +1-n$ and $g=\dim \Gamma(C, \omega_C)$ (the genus of $C$).
Thus we see that the degree of $L_C$ is greater than $2g$. From Theorem 6 in \cite{Mf},
we see that $L_C$ is normally generated.  It is easy to get the normal generation of $L$
from that of $L_C$.  \hfill $\Box$

\bigskip

Next we define the notion of a nonsingular or Gorenstein vertex of a lattice polytope.

Let $M$ be a free abelian group of rank $n\ge2$ and $P$ a lattice $n$-polytope
in $M_{\mathbb{R}}$.  For a vertex $v$ of $P$, make the cone
$$
C_v(P):=\mathbb{R}_{\ge0}(P-v) =\{r(x-v)\in M_{\mathbb{R}}; \ \mbox{$r\ge0$ and $x\in P$}\}.
$$
$P$ is called {\it Gorenstein at $v$} if there exists a lattice point $m_0$ in $C_v(P)$ such that
the equality
\[
(\mbox{Int}C_v(P))\cap M =m_0 +C_v(P)\cap M
\]
holds.   $P$ is called {\it Gorenstein} if it is Gorenstein at all vertices.

Let $M^*:=\mbox{Hom}_{\mathbb{Z}}(M, \mathbb{Z})$ be the dual to $M$ with the natural
pairing $\langle, \rangle: M^*\times M\to \mathbb{Z}$.
We define the dual cone of $C_v(P)$ as
\[
C_v(P)^{\vee}:=\{y\in M^*_{\mathbb{R}}; \langle y, \ \mbox{$x\rangle \ge0$ for all $x\in
C_v(P)$}\}.
\]
Since $C_v(P)^{\vee}$ is a polyhedral cone, there are faces of dimension 1, which are called
{\it rays}, $\rho_1, \dots, \rho_s$ such that 
\[
C_v(P)^{\vee}=\rho_1+\dots +\rho_s.
\]
For a ray $\rho_i$, set $u_i\in M^*$ the generator of the semi-group $\rho_i\cap M$.
We know that the vertex $v$ of $P$ is Gorenstein if and only if there exists a lattice point
$m_0\in M$ such that $\langle u_i, m_0\rangle=1$ for $i=1, \dots, s$.
Moreover, $v$ is
called {\it nonsingular} if
 $s=n$ and if $\{u_1, \dots, u_n\}$ is a $\mathbb{Z}$-basis of $M^*$.

In the following sections, we often use three lemmas proved in \cite{OZ}.
\begin{lem}\label{l1}
Let $P$ be a lattice polytope of dimension $n$.
If there exists an integer $r$ with $1\le r\le n-1$ satisfying the condition that
the multiple $rP$ does not contain lattice points in its interior,
then the equality
$$
(kP)\cap M +P\cap M =((k+1)P)\cap M
$$
holds for all integers $k\ge n-r$.
\end{lem}

If a lattice $n$-simplex is isomorphic to the convex hull of a $\mathbb{Z}$-basis of $M$,
it is called {\it basic}.

\begin{lem}\label{l2}
Let $P$ be a lattice  $n$-simplex with $\sharp\{P\cap M\}=n+1$.
If the multiple $(n-1)P$ does not contain lattice points in its interior,
then it is basic.
\end{lem}

\begin{lem}\label{l3}
Let $P$ be a lattice polytope of dimension $n$.
If the multiple $nP$ does not contain lattice points in its interior,
then it is basic.
\end{lem}

\section{Lattice Polygons}

Let $M$ be a free abelian group of rank two and $P$ a lattice $2$-polytope in
$M_{\mathbb{R}}$, which is often called {\it a lattice polygon}.

\begin{prop}\label{2;p1}
Let $P$ be a lattice polygon such that the number of lattice points contained in its interior
is one.
Then $P$ is Gorenstein.
\end{prop}
{\it Proof.}
We assume that $P\cap M=\{m_0\}$.
Let $v$ be a vertex of $P$ and $E_1, E_2$ two edges starting from $v$.
Let $m_i\in E_i\cap M$ be the lattice point nearest to $v$ for $i=1, 2$.
By taking a suitable affine transformation of $M$, we may set $v=0, 
m_0=(1,1), m_1=(1,0)$
and $m_2=(a, b)$ for $0\le a<b$.

If $a=0$, then $v$ is a nonsingular vertex.

Set $a\ge1$.  We know that if $b=a+1$, then $v$ is a Gorenstein vertex.

If $a=1$, then $b=2$ because $(\mbox{Int}P)\cap M=\{(1,1)\}$.

If $a\ge2$, then $b< 2a$ since $(1,2)\notin P$.
 If $a=2$, then $b=3$ because ${\rm g.c.d}.(a,b)=1$.
 
 If $a\ge3$, then $P$ has another edge $E_3$ starting from $(1,0)$ and
 it does not contain $(2,2)$ in its interior, hence, we have $b\ge 2a-2$.
 We see that for $a\ge3$ $b=2a-2$, or $b=2a-1$.
 
 If $a\ge3$ and $b=2a-1$, then $(2,3)$ is contained in the interior of the triangle 
$\mbox{Conv}\{v, m_1, m_2\}$.

If $(a,b)=(3,4)$, then the triangle  $\mbox{Conv}\{v, m_1, m_2\}$ is Gorenstein.

If $a\ge4$ and $b=2a-2$, then $a$ would be an odd integer, that is, $a\ge5$.
In this case the triangle  $\mbox{Conv}\{v, m_1, m_2\}$ contains $(2,3)$
in its interior.

Thus we see that if $v$ is a singular vertex, then $(a,b)=(1,2), (2,3)$, or $(3,4)$.
In any case $v$ is a Gorenstein vertex. \hfill $\Box$

\bigskip

We note that the same statement in dimension three as in Proposition~\ref{2;p1}
does not hold.  For example, consider the tetrahedron
$D:=\mbox{Conv}\{0, (1,0,0), (0,1,0), (2,2,5)\}$, which contains $(1,1,2)$ in its interior.
The vertex $0$ of $D$ is not a Gorenstein vertex because $(1,1,1)$ is not contained in $D$.

\section{Pyramids and $3$-polytopes}

Let $M$ be a free abelian group of rank $n\ge3$.
A lattice $n$-polytope $P$ in $M_{\mathbb{R}}$ is called {\it a pyramid} if
there exists a facet $F$ and a vertex $v\notin F$ such that $P=\mbox{Conv}\{F, v\}$.
 A typical example of a pyramid is a lattice $n$-simplex, that is, it is the convex hull of
 $(n+1)$ lattice points in general position.
 We also have a Gorenstein lattice $n$-simplex with index $n-1$.
 We define as
 \[
 D_n:=\mbox{Conv}\{0, e_1, e_2, e_1;e_2+2e_3, e_4, \dots, e_n\},
 \]
 for a $\mathbb{Z}$-basis $\{e_1, \dots, e_n\}$ of $M$.
 We note that $D_n$ is Gorenstein and it is not normal.
 
 In \cite{OZ} Ogata and Zhao proved Theorem~\ref{0;t01}, whose statement can
 be interpreted into the following Proposition.
 \begin{prop}[Ogata-Zhao]\label{3;p1}
Let $D$ be a lattice $n$-pyramid with $\sharp(D\cap M)=n+1$ for $n\ge3$.
If $\sharp\{\mbox{\rm Int}((n-1)D)\cap M\}=1$, then $D$ is isomorphic to $D_n$.
\end{prop}

For general pyramids we give a characterization of Gorenstein polytopes.

\begin{lem}\label{3;l2}
Let $P=\mbox{\rm Conv}\{F, v\}$ be a lattice $n$-pyramid for $n\ge3$.
Assume that $P\cap M=(F\cap M)\cup\{v\}$ 
and $\sharp\{\mbox{\rm Int}((n-1)P)\cap M\}=1$.
If $F$ is Gorenstein, then so is  $P$.
\end{lem}
{\it Proof.}
If $\sharp(P\cap M)=n+1$, then it is Gorenstein from Proposition~\ref{3;p1}.

Assume that $\sharp(P\cap M)\ge n+2$.  Then $P$ is normal from Proposition~\ref{1;p1}.
Set $M':=(\mathbb{R}(F-v))\cap M$.  Then $M'$ has rank $n-1$.  Since $P$ is normal,
$M'$ is a direct summand of $M$, that is, $M\cong M'\oplus \mathbb{Z}m$
for some $m\in M$.
If we set as $\{m_0\}=\mbox{\rm Int}((n-1)P)\cap M$, then $m_0$ is contained in
the interior of $(n-2)F$.

Since $F$ is Gorenstein, a vertex $v'$ of $F$ is a Gorenstein vertex of $P$.
Since
\[
C_v(P)\cap M \cong \bigoplus_{t\ge 0}C_v(P)\cap (M'\oplus tm)
\]
and since $\mbox{Int}(tF)\cap M=m_0+((t+2-n)F)\cap M$ for $t\ge n-1$, the vertex $v$ is 
Gorenstein. \hfill $\Box$

\medskip

\begin{cor}\label{3;c1}
Let $P=\mbox{\rm Conv}\{F, v\}$ be a lattice $3$-pyramid with $P\cap M=(F\cap M)\cup\{v\}$.
If $\sharp\{\mbox{\rm Int}(2P)\cap M\}=1$, then $P$ is Gorenstein.
\end{cor}
{\it Proof.}
Since $\sharp\{(\mbox{Int}P)\cap M\}=1$, we see that $F$ is a Gorenstein lattice polygon from
Proposition~\ref{2;p1}. \hfill $\Box$

\medskip

Now we give a characterization of Gorenstein $3$-polytopes.

\begin{prop}\label{3;p1}
Let $P$ be a lattice $3$-polytope in $M_{\mathbb{R}}$.
If $\sharp\{\mbox{\rm Int}(2P)\cap M\}=1$, then $P$ is Gorenstein.
\end{prop}
{\it Proof.}
Let $v$ be a singular vertex of $P$. We may assume that the set of all rays of the cone
$C_v(P)$ is $\{\rho_1, \dots, \rho_t\}$ with $t\ge3$.
Let $m_i$ be the generator of the semi-group $\rho_i\cap M$ for $1\le i\le t$.
Set $Q:=\mbox{Conv}\{0, m_1, \dots, m_t\}$.

We separate our argument into two cases: 

(a) The case when $Q$ has a facet $G$ containing all $m_1, \dots, m_t$, that is,
$Q$ is a pyramid $\mbox{Conv}\{v, G\}$.  We may assume $Q\not\cong D_3$.
Since $\sharp\{\mbox{Int}(2Q)\cap M\}\le1$, we see that $Q$ is normal from Lemma~\ref{l1} and Proposition~\ref{2;p1}.

If $Q\cap M =(G\cap M)\cup\{v\}$, then we have a direct sum decomposition
$M\cong M'\oplus \mathbb{Z}m$ with $G-v\subset M'_{\mathbb{R}}$ as in the proof of
Lemma~\ref{3;l2}.
Thus, if $\sharp\{\mbox{Int}(2Q)\cap M\}=1$, then $\mbox{Int}(2Q)\cap M=
(\mbox{Int}G)\cap M$, hence, $G$ is a Gorenstein polygon.
If $\sharp\{\mbox{Int}(2Q)\cap M\}=0$, then $\sharp\{\mbox{Int}(2G)\cap M\}=1$
since $G$ is not a facet of $P$.
In any case $G$ is a Gorenstein polygon.   Hence, $v$ is a Gorenstein vertex of $Q$.

If $Q\cap M \not=(G\cap M)\cup\{v\}$, we may assume that a facet $F_0:=
\mbox{Conv}\{0, m_1, m_2\}$ contains a lattice point $m'$ in its relative interior.
Since the lattice point $m'+m_i$ is contained in the interior of $2Q$ for $3\le i\le t$, we see
$t=3$.
By the same reason we see that $\sharp\{(\mbox{Int}F_0)\cap M\}=1$,
$\mbox{Conv}\{0, m_i, m_3\}$ contains no lattice points 
in its interior for $i=1, 2$ and that two segments $[m_1, m_3]$ and $[m_2, m_3]$ have lattice length one.  
Since $Q$ is a pyramid of the form $\mbox{Conv}\{m_3, F_0\}$ 
and since $Q\cap M=(F_0\cap M)\cup\{m_3\}$, 
we see that $Q$ is Gorenstein from Lemma~\ref{3;l2}, hence, 
 $v$ is a Gorenstein vertex of $P$.

(b) The case when $Q$ is not a pyramid.
We may assume that the segment $[m_1, m_s]$ is the edge of $Q$ with $3\le s<t$.
Then the relative interior of the segment $[m_2, m_t]$ is contained in the interior of $Q$ and
the lattice point $m_2+m_t$ is contained in the interior of $2Q$.
Thus we see that $s=3$ and $t=4$.
By the same reason, we see that four edges $[m_1, m_i]$ and $[m_3, m_i]$ for $i=2,4$ 
have lattice length one.
Since $\mbox{Int}(2Q)\cap M=\{m_2+m_4\}$, two lattice 3-simplices $Q_i:=\mbox{Conv}
\{0, m_j, m_2, m_4\}$ for $j=1,3$ are basic by Lemma~\ref{l2}. 
 Thus $v$ is a Gorenstein vertex of $Q$.

In both cases (a) and (b), we prove that $v$ is a Gorenstein vertex of $P$.
\hfill $\Box$

\section{Lattice polytopes in higher dimension}\label{sect1}

Let $M$ be a free abelian group of rank $n\ge4$.
Let $P$ be a lattice $n$-polytope with $\sharp (\mbox{\rm Int}(n-1)P)\cap M=1$.
Take a vertex $v$ of $P$.
Let $\{\rho_1, \dots, \rho_t\}$ be the set of all rays of the cone $C_v(P)=\mathbb{R}_{\ge0}
(P-v)$ and let $m_i$ the generator of the semi-group $\rho_i\cap M$ for $i=1, \dots, t$.
Set $Q:=\mbox{Conv}\{0, m_1, \dots, m_t\}$ as in the proof of Proposition~\ref{3;p1}.

\begin{lem}\label{4;l1}
In the above notation, $Q$ is a pyramid.
\end{lem}
{\it Proof}.
If $Q$ is not a pyramid, then there exist two facets $G_1, G_2$ 
with $\dim G_1\cap G_2=n-2$ of $Q$ not containing $v$.
Among $\{m_1, \dots, m_t\}$ choose a $m_i$ in $G_1\setminus G_2$ and $m_j$ in
$G_2\setminus G_1$.  Then the relative interior of the segment $[m_i, m_j]$ is
contained in the interior of $Q$ and, hence the lattice point $m_i+m_j$ is contained in
$\mbox{Int}(2Q)$.  Since $n-2\ge2$ and since $\mbox{Int}((n-2)Q)\cap M=\emptyset$,
it is impossible.  \hfill $\Box$

\medskip
From this Lemma we obtain a characterization of Gorenstein polytopes.

\begin{prop}\label{4;p1}
Let $P$ be a lattice $n$-polytope for $n\ge4$.
If $\sharp\{(\mbox{\rm Int}(n-1)P)\cap M\}=1$, then $P$ is Gorenstein.
\end{prop}
{\it Proof.}
Let $v$ be a vertex of $P$.   Make the cone $C_v(P)$ with the apex $v$ as above.
Let $\{m_1, \dots, m_t\}$ be the primitive minimal generator of $C_v(P)$.
Set $Q=\mbox{Conv}\{0, m_1, \dots, m_t\}$.
We may assume that $v$ is a singular vertex and that $Q$ is not isomorphic to $D_n$.

By Lemma~\ref{4;l1}, $Q$ is a pyramid, that is, there exists a facet $G$ containing all
$m_i$'s such that $Q=\mbox{Conv}\{v, G\}$.  By assumption 
$\sharp\{(\mbox{\rm Int}(n-1)Q)\cap M\}\le1$.
We see that $Q$ is normal from Lemma~\ref{l1} and Proposition~\ref{1;p1}.

(I) The case when $Q\cap M=(G\cap M)\cup\{v\}$.
Then we have a direct sum decomposition $M\cong M'\oplus \mathbb{Z}m$ with
$G-v\subset M'_{\mathbb{R}}$.
Thus, $\sharp\{\mbox{Int}(n-2)G\cap M\}=1$ if $\sharp\{(\mbox{\rm Int}(n-1)Q)\cap M\}=1$.
If $\sharp\{(\mbox{\rm Int}(n-1)Q)\cap M\}=0$, then $\sharp\{\mbox{Int}(n-1)G\cap M\}=1$
because $G$ is not a facet of $P$ and $(\mbox{\rm Int}(nQ))\cap M\not=\emptyset$
by Lemma~\ref{l3}.
If $\sharp\{\mbox{Int}(n-2)G\cap M\}=1$, then $G$ is a Gorenstein $(n-1)$-polytope by
the induction hypothesis on dimension.  If $\sharp\{\mbox{Int}(n-1)G\cap M\}=1$, then
$G$ is also  Gorenstein by Proposition 1 in \cite{OZ}.
From Lemma~\ref{3;l2} we see that $Q$ is Gorenstein.

(II) The case when $Q\cap M\not=(G\cap M)\cup\{v\}$.
We may assume that a face $E=\mbox{Conv}\{0, m_1, \dots, m_s\}$ of $Q$ contains a
lattice point $m'$ in its relative interior by renumbering $m_i$'s.  We note $2\le \dim E\le n-1$. 
Set $r=\dim E$.  
Since $\dim Q=n$, we can choose $n-r$ from $\{m_{s+1}, \dots, m_t\}$ such that
the sum of that with $m'$ is a lattice point in the interior of $(n+1-r)Q$.  Thus $r=2$ and $s=2$.
 By renumbering we may assume that $m'+m_3+\dots +m_n$ is in the interior of $(n-1)Q$.
 If $t>n$, then $m_t$ is not contained in the $(n-1)$-simplex $\mbox{Conv}\{m_1, \dots, m_n\}$.
 If the segment $[m', m_t]$ has an intersection with the interior of $Q$, then  the lattice
 point $m'+m_t$ is contained in $\mbox{Int}(2Q)$, which contradicts with
 $\mbox{Int}(n-2)Q\cap M=\emptyset$.
 If there exists a facet of $Q$ containing $E$ and $m_t$, we may assume that $m_t$
 locates in the opposite side to $m_n$ with respect to the hyperplane
 through $\{0, m_1, \dots, m_{n-1}\}$.  Then $m'+m_3+\dots +m_{n-1}+m_t$ would be
 another lattice point in the interior of $(n-1)Q$.
 This contradicts with $\sharp\{\mbox{Int}(n-1)Q\cap M\}=1$.
 Thus we have $t=n$.
 By the same reason we see that $\sharp\{(\mbox{Int}E)\cap M\}=1$ and that such $E$ is
 unique.
 
 We claim $\sharp\{G\cap M\}=n$.
 If a face $E'$ of $G$ contains a lattice points in its relative interior, then $2\mbox{Conv}
 \{E, E'\}$ contains lattice points more than one in its relative interior.
 This contradicts with $\sharp\{(\mbox{Int}(n-1)Q)\cap M\}=1$.
 
 Since $\mbox{Int}(n-2)G\cap M=\emptyset$, we see that $G$ is a basic $(n-1)$-simplex
 from Lemma~\ref{l2}.
 Set $G':=\mbox{Conv}\{0, m_1, \dots, m_{n-1}\}$.  The lattice point $m'+m_3+\dots +m_{n-1}$
 is a unique lattice point contained in the interior of $(n-2)G'$.  By the induction hypothesis,
 $G'$ is a Gorenstein $(n-1)$-polytope.
 We may write as $Q=\mbox{Conv}\{m_n, G'\}$, which is a pyramid with
 $Q\cap M=(G'\cap M)\cup\{m_n\}$.   By (I), we see that $Q$ is Gorenstein, hence,
 $v$ is a Gorenstein vertex of $P$.
\hfill $\Box$

\end{document}